\documentclass[a4paper,11pt]{article} 
\usepackage{amsmath,amssymb,amsfonts} 
\usepackage[all]{xy}

\newcommand{\CC}{\mathbb{C}}

\newcommand{\PP}{\mathbb{P}}
\newcommand{\QQ}{\mathbb{Q}}
\newcommand{\ZZ}{\mathbb{Z}}
\newcommand{\Ay}{\mathcal{A}}

\newcommand{\Oh}{\mathcal{O}}

\newcommand{\eM}{\mathcal{M}}

\newcommand{\GL}{\mathrm{GL}}
\newcommand{\al}{\alpha}
\newcommand{\be}{\beta}
\newcommand{\gam}{\gamma}
\newcommand{\del}{\delta}

\newcommand{\fie}{\varphi}
\newcommand{\Fie}{\Phi}

\newcommand{\la}{\lambda}

\newcommand{\Ttilde}{\widetilde T}

\newcommand{\Vhat}{\widehat V}
\newcommand{\What}{\widehat W}
\newcommand{\vhat}{\widehat v}
\newcommand{\what}{\widehat w}
\newcommand{\nuhat}{\widehat\nu}

\DeclareMathOperator{\Proj}{Proj}

\DeclareMathOperator{\coker}{coker}

\newtheorem{thm}{Theorem}[section]
\newtheorem{mainthm}[thm]{Main Theorem}

\newtheorem{cor}[thm]{Corollary}

\newtheorem{rmk}[thm]{Remark}
\newenvironment{pf}{\paragraph{Proof}}{\par\medskip}

\newenvironment{pfthm}{\paragraph{Proof of main theorem}}{\par\medskip}

\title{Extending symmetric determinantal\\quartic surfaces}
\author{Stephen Coughlan}
\date{}

\begin{document}
\maketitle
\begin{abstract}We give an explicit construction for the extension of a symmetric determinantal quartic K3 surface to a Fano $6$-fold. Remarkably, the moduli of the $6$-fold extension are in one-to-one correspondence with the moduli of the quartic surface. As a consequence, we determine a $16$-parameter family of surfaces of general type with $p_g=1$ and $K^2=2$ as weighted complete intersections inside Fano $6$-folds.\end{abstract}
\section{Introduction}
Let $D$ be a curve of genus $3$ which is not hyperelliptic. Then the canonical model of $D$ is a plane quartic, and any such quartic has $36$ ineffective theta characteristics, of which we fix one and call it $A$. In this paper we study extensions of the graded ring \[R(D,A)=\bigoplus_{n\ge0} H^0\left(D,\Oh_D(nA)\right),\] where $\Proj R(D,A)$ gives $D\subset\PP(2^3,3^4)$. The structure of $R(D,A)$ is completely determined by a symmetric $4\times4$ matrix with linear entries in $3$ variables $y_1,y_2,y_3$ of weight $2$, hence we call $D$ a \emph{symmetric determinantal quartic} curve. If we add another variable $y_0$ of weight $2$ into the matrix, preserving the linearity and symmetry, we get the graded ring of a symmetric determinantal K3 surface $T\subset\PP(2^4,3^4)$ with $10\times\frac12(1,1)$ points. A priori we know that both $T$ and $D$ are always symmetric determinantal varieties, so this is the only way to extend $D$ to a K3 surface with $10\times\frac12$ points. See Section~\ref{sec!symdet} of this paper for further remarks on symmetric determinantal varieties.

Now, since $T$ is a K3 surface, it is naturally the elephant hyperplane section of a Fano $3$-fold $W\subset\PP(1,2^4,3^4)$ with $10\times\frac 12$ points. In other words, $T$ is the hyperplane section of weighted degree $1$ \[T=H\cap W\subset\PP(1,2^4,3^4),\] or in terms of graded rings, there is an element $a\in H^0(W,\Oh({-}K_W))$ whose vanishing defines $H$ and so \[R(T,A)=R(W,{-}K_W)/(a).\] This process can be iterated and we can continue incorporating more variables $b$, $c$, $d$ of degree $1$ into the ring. We obtain a tower of inclusions \[D\subset T\subset W^3\subset W^4\subset W^5\subset W^6\subset\PP(1^4,2^4,3^4),\] where each $W^n$ is a Fano $n$-fold of Fano index $n-2$. Having built the tower as far as a Fano $6$-fold, we discover an amazing one-to-one correspondence between the moduli of the K3 surface $T$ and the moduli of the $6$-fold $W^6$.
\begin{mainthm}\label{thm!nonhell} For each quasismooth symmetric determinantal K3 surface $T\subset\PP(2^4,3^4)$ with $10\times\frac12$ points there is a unique extension to a quasismooth Fano $6$-fold $W\subset\PP(1^4,2^4,3^4)$ with $10\times\frac12$ orbifold points and such that \[T=W\cap H_1\cap H_2\cap H_3\cap H_4,\] where the $H_i$ are hyperplanes of the projective space $\PP(1^4,2^4,3^4)$.\end{mainthm}

Jan Stevens \cite{S}, first observed this phenomenon in 1993 when calculating the deformation--extension theory for the special case of the Klein quartic curve, which has maximal symmetry group of order $168$. This extra symmetry restricts the deformation extension space enough to make the computation viable. It is not immediately clear how to perform this extension procedure in general; we believe it is not as simple as generalising the symmetric matrix to have entries involving $a,\dots,d$.

We prove the theorem in Section~\ref{sec!symdetextend} for any symmetric determinantal quartic surface. The novel idea is to consider the image of the symmetric determinantal surface $T$ under a projection map, extend this image, then reverse the $6$-dimensional projection to obtain $W$. The advantage of our approach is that we sidestep the complicated calculations involved in extending $T$ directly, the disadvantage is that so far we have not been able to explain the structure of the $6$-fold $W$ in terms of the symmetric matrix.

The Fano $6$-fold $W$ is an example of a \emph{key variety}: lots of interesting varieties are contained in $W$ as appropriate weighted complete intersections. We have already seen how to obtain $T$ from $W$, and the curve $D$ of genus $3$ is obtained as $T\cap Q$ where $Q$ is a hypersurface of weight $2$ avoiding the $\frac12$ points. A simple application of this key variety principle leads us to an important family of surfaces of general type.

\begin{cor} There is a $16$-parameter family of surfaces $Y$ of general type with $p_g=1$, $q=0$, $K^2=2$ and no torsion, each of which is a complete intersection of type $(1,1,1,2)$ in a Fano $6$-fold $W\subset\PP(1^4,2^4,3^4)$ with $10\times\frac12$ points.\end{cor}

The proof of this corollary is in Section~\ref{sec!surfaces}. We observe that the expected dimension of the moduli space of surfaces $Y$ of general type with $p_g=1$, $q=0$ and $K_Y^2=2$ is $16$, and that the unique curve $D\in|K_Y|$ is precisely the symmetric deteminantal curve of genus $3$ described above. Such surfaces $Y$ were constructed by Catanese and Debarre in \cite{CDE}, by examining the image of the bicanonical map as suggested by Enriques and splitting into cases accordingly. Todorov \cite{To} also studied the case with torsion $\ZZ/2$, and more recently examples have been constructed using $\QQ$-Gorenstein smoothing theory in \cite{PPS}. Our method is new and is more widely applicable to other examples. There is also a hyperelliptic degeneration of this construction, which has applications to Godeaux surfaces with $\ZZ/2$-torsion. We return to this topic in the forthcoming paper \cite{part2}.

\subsection*{Acknowledgements} I would like to thank Miles Reid for introducing me to this problem, which forms part of my University of Warwick PhD thesis \cite{mythesis}; in addition I thank Jan Stevens for his useful comments and suggestions. This research was partially supported by the World Class University program through the National Research Foundation of Korea funded by the Ministry of Education, Science and Technology (R33-2008-000-10101-0).

\section{Symmetric determinantal varieties}\label{sec!symdet}
In this paragraph we collect together various facts about symmetric determinantal varieties and ineffective theta characteristics. Of particular importance is the projection construction of Section~\ref{sec!projection}, which is used in the proof of Main Theorem \ref{thm!nonhell}.

Take $D_4\subset\PP^2$ the canonical model of a genus $3$ curve, and let $A$ be an \emph{ineffective theta characteristic} on $D$. That is, a divisor class $A$ such that $2A=K_D$ and $h^0(A)=0$. We know that $A$ exists because there are $28$ bitangents $\be_i$ to $D$, and any of the $36$ combinations $\be_j-\be_k+\be_l$ is an ineffective theta characteristic.

We explain how the existence of $A$ is equivalent to a symmetric determinantal representation of $D$ by using well known results on projectively Cohen--Macaulay sheaves. A coherent sheaf $\Ay$ on $\PP^n$ is called projectively Cohen--Macaulay if its associated module $\Gamma_*\Ay$ is a Cohen--Macaulay graded $k[\PP^n]$-module. Thus if the support of $\Ay$ is a hypersurface $X\subset\PP^n$, then using the free resolution of $\Gamma_*\Ay$ we get a locally free resolution of $\Ay$ \begin{equation}\label{res!A}0\leftarrow\Ay\leftarrow\bigoplus_{i=1}^m\Oh_{\PP^n}(-d_i)\xleftarrow{\,M\,}\bigoplus_{i=1}^m\Oh_{\PP^n}(-e_i)\leftarrow0,\end{equation}
where the vanishing of the determinant of $M$ defines $X$ set-theoretically, and the degree of $X$ is $\sum(e_i-d_i)$. If in addition we require \[H^0(\Ay(-1))=H^{n-1}(\Ay(1-n))=0,\] then all the $d_i=0$, $e_i=1$ so that $M$ has linear entries and $X$ has degree $m$. Finally, $M$ will be symmetric if $\Ay(-1)^{[2]}=\Oh_X(1)$. These conditions are clearly satisfied in the case $X$ is a plane quartic, and $\Ay=\Oh_D(A)(1)$. Conversely if $D_4\subset\PP^2$ is defined by the determinant of a $4\times4$ symmetric matrix $M$ with linear entries then $\Ay:=\coker M$ is a projectively Cohen--Macaulay sheaf on $D_4$ with $\Ay(-1)^{[2]}=\Oh_D(1)$. See \cite{Be} for details on projectively Cohen--Macaulay sheaves in this context.

Now writing $\Oh_D(A)=\Ay(-1)$ and twisting by $-1$ the short exact sequence (\ref{res!A}) becomes \[0\leftarrow\Oh_D(A)\xleftarrow{(z_i)}4\Oh_{\PP^2}(-1)\xleftarrow{M}4\Oh_{\PP^2}(-2)\xleftarrow{(z_i)^t}0.\] From this we deduce that the curve \[D=\Proj R(D_4,A)\subset\PP(2^3,3^4)\] has equations \begin{equation}\label{eq!symminors}\left(\begin{matrix}z_1,&z_2,&z_3,&z_4\end{matrix}\right)M=0,\quad\bigwedge^3_{i,j}M=z_iz_j,\end{equation} where $\bigwedge^3_{i,j}M=(-1)^{i+j}\det M_{ij}$, the $(i,j)$th cofactor of $M$. See \cite{Cat1} for a proof of this.

In a similar manner, if we allow $M$ to be symmetric but with linear entries in four variables $y_1,\dots,y_4$ then the projective variety $T_4$ defined by $\det M=0$ is a quartic K3 surface in $\PP^3$ with $10$ nodes. The above properties of projectively Cohen--Macaulay sheaves imply that there is an ineffective divisor class $A$ on $T$ such that $\Oh_T(A)^{[2]}=\Oh_T(1)$. There is a similar short exact sequence and the equations of $T\subset\PP(2^4,3^4)$ are also the same, but now $M$ has entries in four variables. The $10$ nodes of $T_4$ become $\frac 12(1,1)$ points of the weighted ambient space.

\subsection{Webs of quadrics}

Let $\PP^9=\PP H^0(\PP^3,\Oh(2))$ be the space of quadrics in $\PP^3$, or if you prefer, the space of symmetric $4\times4$ matrices up to scalar multiplication. There is a natural stratification of this space by rank: \[\PP^9\supset V^8_4\supset V^6_{10}\supset V^3_8.\] For example, $V^8_4$ is a hypersurface of degree $4$ in $\PP^9$, which corresponds to quadrics in $\PP^3$ of rank $\le 3$, or equivalently $4\times4$ symmetric matrices whose determinant vanishes. Similarly $V^6_{10}$ (respectively $V^3_8$) is the locus of quadrics of rank $\le2$ (resp. $\le1$).

Now take a web $\eM$ of quadrics in $\PP^3$, i.e.~$\eM$ is a linear system of projective dimension $3$ inside $\PP H^0(\PP^3,\Oh(2))$. Choose coordinates $y_1,\dots,y_4$ for $\eM$, and define \[T_4=\eM\cap V^8_4.\] Then $T_4$ is the locus of quadrics of rank $\le3$ in $\eM$, and it is defined by the vanishing of the determinant of a $4\times4$ symmetric matrix with linear entries in $y_1,\dots,y_4$. In general, if $\eM$ is base point free then this is an irreducible quartic hypersurface in $\PP^3$.

The singularities of $T_4$ are given by $T_4\cap V^6_{10}$, the locus of quadrics of rank $\le2$ in $\eM$. There are $10$ isolated points in this locus, corresponding to $10$ nodes on $T_4$. Of course $T_4\cap V^6_{10}$ is defined algebraically by the vanishing of the $3\times3$ minors of $M$. These minors generate the linear system of contact cubics to the quartic hypersurface $T_4$. See \cite{Tj} for details of this and \cite{Ca} for a classical proof.

\subsection{An almost homogeneous space}
Let $V=\CC^4$ be a vector space of dimension $4$, then there is a natural $G=\GL(4,\CC)$ group action on $V$ by matrix multiplication, and this induces an action of $G$ on the vector spaces $S^2(V)$ and $\bigwedge^3V$. We define the almost homogeneous space $X$ to be the closure of the $G$-orbit of the vector
\[\left(A,B\right)\in S^2(V)\oplus\bigwedge^3V\] where \[A=\begin{pmatrix}1&&&\\&1&&\\&&1&\\&&&0\end{pmatrix}, \quad B=\begin{pmatrix}0\\0\\0\\1\end{pmatrix}.\] The K3 surface $T$ is the intersection of $X$ with a $4$-dimensional subspace $\eM\subset S^2(V)$ and naturally lives in weighted projective space $\PP(2^4,3^4)$ with equations (\ref{eq!symminors}).

\subsection{A projection construction for $T$}\label{sec!projection} Let $T_4\subset\PP^3$ be a hypersurface with $10$ nodes. Choose a node and project away from it onto the complementary plane $\PP^2$. Explicitly, we can choose coordinates so that the equation of $T_4$ is \[\al_2(y_1,y_2,y_3)y_4^2+\be_3(y_1,y_2,y_3)y_4+\gam_4(y_1,y_2,y_3)=0,\] with a node at $P=(0,0,0,1)$. Then linear projection onto the plane with coordinates $y_1$, $y_2$, $y_3$ gives a double covering of the plane branched in the sextic curve $\be^2-4\al\gam$. The image of $P$ under the projection is the conic $\al=0$, which touches the branch curve doubly in each of $6$ points. We say that the conic is \emph{totally tangent} to the sextic.

If we further assume that $T$ is a symmetric determinantal hypersurface, then an explicit calculation shows that the branch curve breaks up into two distinct cubics. These two cubics intersect one another transversally to give $9$ nodes, and the additional node from the centre of projection makes $10$ nodes on $T_4$.

The same map can also be viewed as a calculation in quasi-Gorenstein projection--unprojection, see \cite{Ki}, \cite{PR} for discussion and examples. Start with the K3 surface $T\subset\PP(2^4,3^4)$ which has $10\times\frac12$ orbifold points. Let $A$ denote the polarising divisor for this model of $T$, choose one of the $\frac12$ points and call it $P$. Then write $\sigma\colon\Ttilde\to T$ for the $(1,1)$-weighted Kawamata blowup of $P$. The exceptional locus $E\cong\PP^1\subset\Ttilde$ is the centre for our projection, and the projection map is determined by the linear system $\sigma^*A-\frac12E$ on $\Ttilde$. The image of this projection is $T^\prime_{6,6}\subset\PP(2^3,3^2)$.

The surface $T^\prime_{6,6}\subset\PP(2^3,3^2)$ is a double cover of $\PP(2,2,2)$ branched in the two cubics defined by the relations of weight $6$. The image of the exceptional curve $E$ is embedded as a conic which is totally tangent to the branch sextic.

A further way to calculate this projection is via explicit commutative algebra. Fairly generally we can assume the matrix $M$ is of the form
\[M=\begin{pmatrix}b&y_4&B&0\\
&a&0&A\\
\text{sym}&&y_1&y_2\\
&&&y_3\end{pmatrix},\]
where $a$, $b$ are general linear forms in $y_1,\dots,y_3$ and $A=y_1+\al_1y_2+\al_2y_3$, $B=\be_1y_1+\be_2y_2+y_3$. The K3 surface determined by this matrix has a $\frac12$ point at $P=(0,0,0,1)$, with local coordinates near the singularity given by the variables $z_3$, $z_4$. Thus if we project away from $P$ we expect to eliminate $y_4$, $z_3$, $z_4$, (c.f.~\cite{Ki}, example~9.13). Calculating cofactors $(1,1)$, $(2,2)$ of $M$ we obtain equations:
\begin{equation}\begin{split}\label{eqn!projected}z_1^2&=F_6(y_1,y_2,y_3)=a(y_1y_3-y_2^2)-y_1A^2,\\
z_2^2&=G_6(y_1,y_2,y_3)=b(y_1y_3-y_2^2)-y_3B^2.
\end{split}\end{equation}
These are the only equations remaining from (\ref{eq!symminors}) that do not involve $y_4$, $z_3$, $z_4$. In particular the cofactor $-M_{12}$ for $z_1z_2$ involves $y_4$ and so does not survive the projection. Further, the product $F_6G_6$ is the defining equation $\be^2-4\al\gam$ of the totally tangent sextic. The equations (\ref{eqn!projected}) define the image of our projection map \[T^\prime_{6,6}\subset\PP(2,2,2,3,3).\]
\begin{rmk}\rm The truncated graded ring $R(T^\prime)^{[2]}$, which is the even subring of $R(T^\prime)$, no longer contains $z_1$, $z_2$ as generators because they have odd degree. However, we win the new generator $z_1z_2$ since we observed above that there is no equation eliminating this product in $R(T^\prime)$. Thus the truncation defines the familiar double cover \[T_{6}\subset\PP(1,1,1,3)\] with equation $u^2=FG$, where $u=z_1z_2$ and we have divided degrees by $2$. Conversely, given $T_6\subset\PP(1,1,1,3)$ defined by $z^2=F_3G_3(y_1,y_2,y_3)$, the above argument shows there is a divisor class $A$ on $T_6$ with $\Oh_{T_6}(2A)=\Oh(1)$.
\end{rmk}

\section{Extending determinantal formats}\label{sec!symdetextend}
In this section we treat extensions of symmetric determinantal quartic surfaces, culminating in the proof of Main Theorem \ref{thm!nonhell}. We use the projection construction for the K3 surface $T$, which is not as symmetric as the determinantal representation but is very beautiful in its own way. Unfortunately a consequence of this approach is that we are not able to completely understand how the symmetric matrix is involved in the extension.

The projection is described in Section~\ref{sec!projection}: we project from one of the $\frac 12$ points on $T$ to get the surface\[\PP^1\xrightarrow{\,\fie\,} T^\prime_{6,6} \subset\PP(2,2,2,3,3)\] with $9\times\frac 12$ points. The surface $T^\prime$ is a double cover of $\PP^2$ branched in a sextic curve which breaks into two cubics. The image of $\fie$ is a conic in the plane $\PP(2,2,2)$ which touches both branch cubics at exactly $3$ points each. Hence constructing a K3 surface $T\subset\PP(2^4,3^4)$ with $10\times\frac12$ points is equivalent to exhibiting a suitable projected surface $T^\prime\subset\PP(2^3,3^2)$ along with a map $\fie$ embedding $\PP^1$ inside $T^\prime$ with appropriate tangency.

Write $y_i$, $z_i$ for the coordinates on $\PP(2,2,2,3,3)$ of weight $2$, $3$ respectively. After coordinate changes, for general $T^\prime$ the embedding of $\PP^1$ is \[\fie\colon\PP^1\to\PP(2,2,2,3,3)\]\begin{equation}\label{eqn!fiemap}(u,v)\mapsto(u^2,uv,v^2,u^3+\al_1u^2v+\al_2uv^2,\be_1u^2v+\be_2uv^2+v^3).\end{equation} We have assumed that $u$ is a factor of $\fie^*(z_1)$ and likewise $v$ divides $\fie^*(z_2)$. Moreover we assume that $\fie^*(z_1)$ and $\fie^*(z_2)$ have no common factor.

Since $S^3(u^2,uv,v^2)$ generates $S^6(u,v)$ we see that the image of $\fie$ is given by the equations
\begin{align}
C_1&\colon z_1^2=y_1(y_1+\al_1y_2+\al_2y_3)^2,\label{eqn!fieimage1}\\
C_2&\colon z_1z_2=y_2(y_1+\al_1y_2+\al_2y_3)(\be_1y_1+\be_2y_2+y_3),\label{eqn!fieimage2}\\
C_3&\colon z_2^2=y_3(\be_1y_1+\be_2y_2+y_3)^2,\label{eqn!fieimage3}\\
Q&\colon y_1y_3=y_2^2.\label{eqn!conic}\end{align} Note that the choice of representation for the first three equations is only unique modulo the conic $Q$ of equation (\ref{eqn!conic}); for example we could have written $z_2^2=\be_1^2y_1^2y_3+2\be_1\be_2y_2^3+(\be_2^2+2\be_1)y_2^2y_3+y_3^3$ instead.

The projected surface $T^\prime$ is given by taking two combinations
\begin{equation}\label{eqn!Tprime}\begin{split}C_1&+l_1(y_1,y_2,y_3)Q\\C_3&+l_3(y_1,y_2,y_3)Q,\end{split}\end{equation} where $l_i$ are linear. There are $9$ moduli for this construction: $3$ from the parameters $\al_i$, $\be_i$ and a further $3$ for each of the linear forms $l_1$, $l_3$. As an illustration, we could choose
\begin{align*}z_1^2&=y_1(y_1+\al_1y_2+\al_2y_3)^2-(y_2+2y_3)(y_1y_3-y_2^2)\\
z_2^2&=y_3(\be_1y_1+\be_2y_2+y_3)^2-y_1(y_1y_3-y_2^2),
\end{align*}
which corresponds to the symmetric matrix
\[M=\renewcommand{\arraystretch}{1.2}\begin{pmatrix}y_1&y_4&\be_1y_1+\be_2y_2+y_3&0\\
&y_2+2y_3&0&y_1+\al_1y_2+\al_2y_3\\
\text{sym}&&y_1&y_2\\
&&&y_3\end{pmatrix}.\]
\begin{rmk}\rm We have made a trade off here between simplifying the equations of $T^\prime$ and simplifying the map $\fie$. Denote the branch cubics by $B_1$, $B_2$ and the conic by $Q$. Then the restrictions $B_i|_Q$ generate a pencil of cubics on $Q\cong\PP^1$. We have chosen $\fie^*(z_i):=B_i|_Q$, which means that the equations of $T^\prime$ take the simpler form $z_i^2=f_i(y_1,y_2,y_3)$. We could have reduced the number of terms involved in the definition of $\fie$ by choosing $\fie^*(z_i)$ to be generators for the pencil of the form $u^3+\al u^2v$ and $\be uv^2+v^3$. However, were we to do this, the price we pay is that we are only able to assume the equations for $T^\prime$ are of the form $(\la_iz_1+\mu_iz_2)^2=f_i(y_1,y_2,y_3)$.
\end{rmk}
\begin{pfthm}
The key point is that there is an analogous projection of the Fano $6$-fold $W$, which has image $\PP^5\subset W^\prime_{6,6}\subset\PP(1^4,2^3,3^2)$. If we can write down the extension of $T^\prime$ to $W^\prime$, then this is as good as extending $T$ to $W$ itself. Of course we have reduced to a much easier problem because we can work explicitly with $T^\prime$ and $W^\prime$ as they are codimension $2$ complete intersections.

We define $\fie$ as in (\ref{eqn!fiemap}) and write $\fie_0\colon\PP^1\to\PP(2,2,2)$ for the standard parametrisation of the conic in $\PP(2,2,2)$: \[\fie_0^*(y_1)=u^2,\quad\fie_0^*(y_2)=uv,\quad\fie_0^*(y_3)=v^2.\] If we write $u$, $v$, $a$, $b$, $c$, $d$ for the coordinates of $\PP^5$ then up to automorphisms of $\PP^5$ and $\PP(1^4,2^3)$, the general extension of $\fie_0$ to $\Fie_0\colon\PP^5\to\PP(1^4,2^3)$ is \[\Fie_0^*(a)=a,\quad\Fie_0^*(b)=b,\quad\Fie_0^*(c)=c,\quad\Fie_0^*(d)=d,\]
\begin{equation}\label{def!Fie0}\setlength\arraycolsep{2pt}\renewcommand{\arraystretch}{1.3}\begin{array}{rllc}\Fie_0^*(y_1)&=u^2&-dv&+bd-c^2,\\\Fie_0^*(y_2)&=uv+bu&+cv&-ad+bc,\\\Fie_0^*(y_3)&=v^2-au&&+ac-b^2.\end{array}\end{equation} The curious extra terms $bd-c^2$, $-ad+bc$, $ac-b^2$ are harmless but they ensure that $\Fie_0^*(y_i)$ are the $2\times2$ cofactors of the matrix \[\begin{pmatrix}a&b-v&c+u\\b+v&c-u&d\end{pmatrix}\] so that the matrix (\ref{matrix!A}) below is more beautiful.

We prove that there is a unique map $\Fie\colon\PP^5\to\PP(1^4,2^3,3^2)$ extending $T^\prime_{6,6}$ to $W^\prime_{6,6}$ and lifting $\Fie_0$ so that the following diagram commutes:
\[\xymatrix@R=30pt@C=35pt
{&\PP(1^4,2^3,3^2)=\Proj S{\ar@{-->}@<-4ex>[d]^{\pi}}\\
\PP^5\ar[ur]+L^<<<<<<<<<{\Fie}\ar[r]^>>>>>>>{\Fie_0}&\PP(1^4,2^3)=\Proj R}\]

Write $M$, $R$, $S$ for the coordinate rings of $\PP^5$, $\PP(1^4,2^3)$ and $\PP(1^4,2^3,3^2)$ respectively. Then $M$ is a graded $R$-module via $\Fie_0^*$ generated by $1$, $u$, $v$ (see equation (\ref{def!Fie0})) with presentation
\[0\leftarrow M\xleftarrow{(1,u,v)} R\oplus 2R(-1)\xleftarrow{A} 2R(-3)\oplus R(-4)\] where $A$ is the matrix 
\begin{equation}\label{matrix!A}\renewcommand{\arraystretch}{1.2}\begin{pmatrix}L_1&L_2&L_3\\
-y_2&y_3&L_2\\
y_1&-y_2&L_1\end{pmatrix}\end{equation}
and the outsized entries are 
\begin{gather*}L_1=by_1+cy_2+dy_3\\L_2=ay_1+by_2+cy_3\\L_3=y_1y_3-y_2^2+b^2y_1+(2bc-ad)y_2+c^2y_3.\end{gather*} Moreover, $M$ is also a graded module over $S$ via $\Fie^*$, with the same generators and of course more relations. Finally, $S$ is a module over $R$ which is not finite. We will not insist on writing $\fie^*$, $\Fie^*_0$, $\Fie^*$ when it is clear that we are dealing with the module structure. 

Since $\Fie$ is a lift of $\Fie_0$ and $\fie$, we can assume the general forms for $\Fie^*(z_i)$ are
\begin{align*}\Fie^*(z_1)&=u^3+\al_1u^2v+\al_2uv^2+s_1u^2+s_2uv+s_3v^2+s_4u+s_5v,\\
              \Fie^*(z_2)&=\be_1u^2v+\be_2uv^2+v^3+t_1u^2+t_2uv+t_3v^2+t_4u+t_5v\end{align*}
where $s_i(a,b,c,d)$, $t_i(a,b,c,d)$ are homogeneous polynomials of degree $1$ or $2$ as appropriate. Now using the $R$-module structure of $M$, we can write
\begin{equation}\begin{split}\label{def!Fie}
\Fie^*(z_1)&=(f+s_4)u+s_5v,\\
\Fie^*(z_2)&=t_4u+(g+t_5)v\end{split}\end{equation} where \[f=y_1+\al_1y_2+\al_2y_3,\quad g=\be_1y_1+\be_2y_2+y_3.\] Here we use coordinate changes such as $z_1\mapsto z_1+s_1y_1$ so that $z_1$, $z_2$ absorb the values of $s_i$, $t_i$ for $i=1$, $2$, $3$. We are required to find suitable values of $s_4$, $s_5$, $t_4$, $t_5$ so that the kernel of $\Fie^*$ contains equations extending (\ref{eqn!fieimage1}), (\ref{eqn!fieimage3}) and (\ref{eqn!conic}). Constructing the extension $\Fie$ of $\fie$ amounts to the following algebraic result:
\begin{thm}\label{thm!extend}
(I)\quad The kernel of $\Fie^*\colon S\to M$ contains equations extending (\ref{eqn!fieimage1}), (\ref{eqn!fieimage3}) of the form 
\begin{align*}z_1^2&-y_1f^2\in R+Rz_1+Rz_2,\\
z_2^2&-y_3g^2\in R+Rz_1+Rz_2
\end{align*}
if and only if $s_4=s_5=t_4=t_5=0$.
\vspace{0.4cm}

\noindent(II)\quad Given part $(I)$, the equations are
\begin{equation}\label{eqn!kerFie1}
z_1^2-y_1f^2=(c^2-bd)f^2-(\del_1L_1-\del_2L_2)df+(\del_1y_2+\del_2y_3)dz_1+\al_2dfz_2\end{equation}
\begin{equation}\label{eqn!kerFie2}
z_2^2-y_3g^2=(b^2-ac)g^2-(-\del_3L_1+\del_1L_2)ag+\be_1agz_1+(\del_3y_1+\del_1y_2)az_2,\end{equation}
where $\del_1=1-\al_2\be_1$, $\del_2=\al_1-\al_2\be_2$, $\del_3=\be_2-\al_1\be_1$.
\end{thm}
\begin{cor}\label{cor!extend} The kernel of $\Fie^*$ contains the following equation extending (\ref{eqn!fieimage2})
\[z_1z_2-fgy_2=fg(ad-bc)-bgz_1-cfz_2,\]
and (nontrivial) equations extending multiples of (\ref{eqn!conic}), of the form
\[y_iL_3\in R+Rz_1+Rz_2\] for $i=1$, $2$, $3$.
\end{cor}
\begin{rmk}\label{rem!unique}\rm Part (I) of the theorem uniquely determines $\Fie$ up to auto\-morphism. Moreover, the coordinate changes used do not alter the original setup \[\fie\colon\PP^1\hookrightarrow T^\prime_{6,6}\subset\PP(2,2,2,3,3),\] so $\Fie$ is completely determined by $\fie$.

As an aside, observe that since we expect the image of $\Fie$ not to be Cohen--Macaulay, the standard strategy of using the hyperplane section principle goes awry. The equation $y_1y_3-y_2^2$ does not extend directly, and we need three separate equations replacing it in the kernel of $\Fie^*$. The image of $\Fie_0\colon\PP^5\to\PP(1^4,2^3)$ is defined by the vanishing of the determinant of the matrix $A$ from (\ref{matrix!A}), which is of degree $8$.
\end{rmk}
\begin{pf} The ``if'' part is a straightforward verification that when $s_4=s_5=t_4=t_5=0$, equations (\ref{eqn!kerFie1}), (\ref{eqn!kerFie2}) are in the kernel of $\Fie^*$ by direct sub\-stitution. The remainder of the proof is the ``only if'' part.

The ring $k[u,v]$ is a graded module over $k[y_1,y_2,y_3]$ via $\fie_0^*$, so referring to equation (\ref{eqn!fiemap}), we can write $\fie^*(z_i)$ as:
\begin{align*}\fie^*(z_1)&=(y_1+\al_1y_2+\al_2y_3)u\\
\fie^*(z_2)&=(\be_1y_1+\be_2y_2+y_3)v.\end{align*}
If we square these two expressions and use the module structure to render residual terms $u^2$, $v^2$ as $y_1$, $y_3$ we obtain the two equations (\ref{eqn!fieimage1}), (\ref{eqn!fieimage3}). Moreover we can write down the equation for $z_1z_2$ by rendering $uv$ as $y_2$.

We attempt the same elimination calculation with $\Fie^*$. Observe that by definition of $\Fie^*$, we can write $u^2$, $uv$, $v^2$ as
\begin{align*}u^2&=\Fie^*(y_1-bd+c^2)+dv\\uv&=\Fie^*(y_2+ad-bc)-bu-cv\\v^2&=\Fie^*(y_3-ac+b^2)+au.\end{align*}
Thus by squaring $\Fie^*(z_i)$ defined in (\ref{def!Fie}) and rendering $u^2$, $uv$, $v^2$ as above, we arrive at
\begin{align*}
\Fie^*\left(z_1^2-f_1^2(y_1-bd+c^2)-2f_1s_5(y_2+ad-bc)-s_5^2(y_3-ac+b^2)\right)&\equiv0\\
\Fie^*\left(z_2^2-t_4^2(y_1-bd+c^2)-2g_1t_4(y_2+ad-bc)-g_1^2(y_3-ac+b^2)\right)&\equiv0
\end{align*} modulo $(a,b,c,d)M$, where $f_1=f+s_4$ and $g_1=g+t_5$. The residual parts to these congruences are
\begin{equation}\label{eqn!res5}\begin{split}K&\colon(f+s_4)^2dv-2(f+s_4)s_5(bu+cv)+s_5^2au\\
L&\colon t_4^2dv-2(g+t_5)t_4(bu+cv)+(g+t_5)^2au,\end{split}\end{equation}
which are homogeneous expressions of degree $6$ in $(a,b,c,d)M$. We prove that for the unique values $s_4=s_5=t_4=t_5=0$ the two residual terms $K$, $L$ are contained in the submodule \[R+Rz_1+Rz_2\subset M=R+Ru+Rv.\] This is necessary and sufficient to obtain equations for $z_i^2$ of the required form in the kernel of $\Fie^*$.

By referring to the definition of $\Fie^*(z_i)$ from (\ref{def!Fie}), we see that the submodule $R+Rz_1+Rz_2$ is the image of the composite map
\[M\xleftarrow{(1,u,v)}R\oplus 2R(-1)\xleftarrow{B} R\oplus 4R(-3)\oplus R(-4)\]
 where $B$ is the matrix \[\renewcommand{\arraystretch}{1.2}
\left(\begin{array}{ccc|ccc}1&0&0&L_1&L_2&L_3\\
                            0&f+s_4&t_4&-y_2&y_3&L_2\\
                            0&s_5&g+t_5&y_1&-y_2&L_1\end{array}\right)\]
Note that the first $3$ columns of $B$ represent the submodule generators $1$, $z_1$, $z_2$ respectively while the last $3$ columns are the matrix $A$ of (\ref{matrix!A}), which maps to $0$ under the composite.

We must write $K$, $L$ of (\ref{eqn!res5}) as expressions in the image of this composite map. We stratify the problem according to degree in $y_1,y_2,y_3$, so that \begin{align*}K&=K^{(2)}+K^{(\le1)},\\L&=L^{(2)}+L^{(\le1)},\end{align*} where for example $K^{(2)}=df^2v$, $L^{(2)}=ag^2u$ are the terms of $K$, $L$ which are degree $2$ in $y_i$. We have to find some $\eta=\eta^{(1)}+\eta^{(0)}$ in $R\oplus 4R(-3)\oplus R(-4)$ such that 
\[K=\begin{pmatrix}1,&u,&v\end{pmatrix}B\eta,\] where $\eta^{(i)}$ has degree $i$ in $y_1,y_2,y_3$.

We can do this explicitly: first work in degree $2$ so that we can assume that the matrix $B$ does not involve $s_i$, $t_i$. We demonstrate how to calculate the preimage $\eta^{(1)}$ of $K^{(2)}=df^2v$ under $B$, as $L^{(2)}$ is exactly similar. The first column of $B$ can be used to eliminate any terms in the first row, so the important part of $B$ is the submatrix \[B^\prime=\begin{pmatrix}f&0&-y_2&y_3\\0&g&y_1&-y_2\end{pmatrix}.\] Since the bottom row of $B^\prime$ only involves $y_3$ as part of $g$, we must write
\begin{align*}f&=y_1+\al_1y_2+\al_2y_3\\
&=y_1+\al_1y_2+\al_2(g-\be_1y_1-\be_2y_2)\end{align*} or as an expression in the bottom row of $B^\prime$,
\[f=\begin{pmatrix}0,&g,&y_1,&-y_2\end{pmatrix}\renewcommand{\arraystretch}{1.2}\begin{pmatrix}*\\\al_2\\1 - \al_2\be_1\\-\al_1+\al_2\be_2\end{pmatrix}.\] We are still free to use the first column of $B^\prime$ to remove spurious terms from the middle row by adjusting the starred entry to solve\[0=\begin{pmatrix}f,&0,&-y_2,&y_3\end{pmatrix}\renewcommand{\arraystretch}{1.2}\begin{pmatrix}*\\\al_2\\1 - \al_2\be_1\\-\al_1+\al_2\be_2\end{pmatrix}.\]This is where we use the extra factor of $f$ in $K^{(2)}$ to avoid having to divide by $f$, so we must have
\begin{align*}\eta_2^{(1)}&=\frac1f(y_2\eta_4^{(1)}-y_3\eta_5^{(1)})&\eta_4^{(1)}&=(1-\al_2\be_1)df\\
\eta_3^{(1)}&=\al_2df&\eta_5^{(1)}&=(-\al_1+\al_2\be_2)df,\end{align*} where $\eta_2^{(1)}$ is the starred entry whose value is completely determined by the rest of $\eta^{(1)}$.
Finally, referring back to the large matrix $B$ and in the same manner as for $B^\prime$, we use the first column to remove any accidental terms from the top row so that the remaining entries of the vector $\eta^{(1)}$ are\begin{align*}\eta_1^{(1)}& =-L_1\eta_4^{(1)}-L_2\eta_5^{(1)}\\\eta_6^{(1)}&=0.\end{align*}

An exactly similar argument proves that \[ag^2u=\begin{pmatrix}1,&u,&v\end{pmatrix}B\xi^{(1)}\] where $\xi^{(1)}$ is the vector \begin{align*}\xi_1^{(1)}&=-L_1\xi_4^{(1)} - L_2\xi_5^{(1)}&\xi_4^{(1)}&=(\be_1\al_1 - \be_2)ag\\
\xi_2^{(1)}&=\be_1ag&\xi_5^{(1)}&=(1 - \be_1\al_2)ag\\
\xi_3^{(1)}&=\frac1g(-y_1\xi_4^{(1)}+y_2\xi_5^{(1)})&\xi_6^{(1)}&=0.\end{align*}

Now we reinstate $s_i$, $t_i$ to the matrix $B$ and use the degree $1$ solutions $\eta^{(1)}$, $\xi^{(1)}$ to compute the full vectors $\eta$, $\xi$. The easiest way to do this is via a direct computation. Evaluate the remaining residual terms 
\begin{align*}K^\prime&:=K-\begin{pmatrix}1,&u,&v\end{pmatrix}B\eta^{(1)}\\
L^\prime&:=L-\begin{pmatrix}1,&u,&v\end{pmatrix}B\xi^{(1)}\end{align*} to obtain two expressions in $M$ of degree $6$ and involving $u$, $v$ in degrees $\le3$. In particular all terms involve some $s_i$ or $t_i$ by construction, and the terms of degree $3$ in $u,v$ have coefficients which must be linear in $s_i$, $t_i$. We attempt to write $K^\prime$, $L^\prime$ as expressions in $R+Rz_1+Rz_2$, first using $z_1,z_2$ to remove terms involving $u^3,v^3$ respectively to obtain $K^{\prime\prime},L^{\prime\prime}:$
\begin{align*}
K^{\prime\prime}&=K^\prime-(-2bs_5 - \al_2dt_4)z_1-\left(2\al_2(-cs_5+ds_4)-\del_2ds_5-\al_2^2dt_5\right)z_2\\
L^{\prime\prime}&=L^\prime-\left(-\be_1^2as_4 - \del_3at_4 + 2\be_1(at_5-bt_4)\right)z_1-(-\be_1as_5-2ct_4)z_2,
\end{align*}
where the $\delta_i$ are the three cross ratios of the $6$ points of tangency on the conic, and they appear in equations~(\ref{eqn!kerFie1}), (\ref{eqn!kerFie2}). Now, in order that $K^{\prime\prime}$ and $L^{\prime\prime}$ are in the submodule $R.1$, the coefficients of $u^2v,uv^2$ occurring in $K^{\prime\prime},L^{\prime\prime}$ must vanish. We write these four coefficients as simultaneous linear equations in the $s_i$, $t_i$
\begin{equation}\label{eqn!siti}C\begin{pmatrix}s_4\\s_5\\t_4\\t_5\end{pmatrix}=0,\end{equation} where $C$ is the coefficient matrix
\[\begin{pmatrix}
\delta_1d& -2\delta_1c + \be_1\delta_2d& 0& -\al_2\delta_1d\\
\delta_2d& -2\delta_2c + (\be_2\delta_2 -\delta_1)d& 0& -\al_2\delta_2d\\
-\be_1\delta_3a& 0& (\al_1\delta_3 -\delta_1)a - 2\delta_3b& \delta_3a\\
-\be_1\delta_1a& 0& \al_2\delta_3a - 2\delta_1b& \delta_1a
\end{pmatrix}.\]
Assume $\Delta$, $\delta_1$ are nonzero\footnote{If $\delta_1=0$ the solution is still $s_i=t_i=0$ but there is an interesting anomaly. See Section~\ref{sec!genposition}.}, where $\Delta=\delta_1^2-\delta_2\delta_3$ is the determinant of the resultant matrix of $f$, $g$ displayed as (\ref{matrix!resultant}) below. Then the unique solution to (\ref{eqn!siti}) is $s_4=s_5=t_4=t_5=0$. Hence $K^{\prime\prime}=L^{\prime\prime}=0$ and so we have proved that $\eta=\eta^{(1)}$ and $\xi=\xi^{(1)}$.

The full form of equation $z_1^2-y_1f^2\in R+Rz_1+Rz_2$ is obtained by writing out the vector $\eta$ inside $R+Rz_1+Rz_2$ in terms of the generators $1$, $z_1$, $z_2$:
\begin{equation*}z_1^2=f^2(y_1-bd+c^2)+\eta_1+\eta_2z_1+\eta_3z_2.\end{equation*} Likewise using $\xi$, the equation for $z_2^2$ is
\begin{equation*}z_2^2=g^2(y_3-ac+b^2)+\xi_1+\xi_2z_1+\xi_3z_2.\end{equation*}
Written out in full, these are equations (\ref{eqn!kerFie1}), (\ref{eqn!kerFie2}) in the statement of the theorem. This concludes the proof of Theorem~\ref{thm!extend}, and its corollary is proved in Section~\ref{sec!pfcor!extend}.\end{pf}

Given the existence of equations extending (\ref{eqn!fieimage1}--\ref{eqn!conic}), we can prove the Main Theorem~\ref{thm!nonhell}: define the unique Fano $6$-fold \[\PP^5\xrightarrow{\Fie}W^\prime_{6,6}\subset\PP(1^4,2^3,3^2)\] extending $\PP^1\xrightarrow{\fie}T^\prime_{6,6}\subset\PP(2,2,2,3,3)$ by taking the combination of equations constructed in Theorem~\ref{thm!extend} and its corollary which correspond to the choice (\ref{eqn!Tprime}) made in the definition of $T^\prime_{6,6}$.
\end{pfthm}
\subsection{General position of tangency points}\label{sec!genposition}
First, if $\Delta=0$ then $\fie^*(z_i)$ have a shared root, which implies one of the tangency points $P$ is common to both branch curves. Thus $P$ is a $\frac12$ point of $T^\prime_{6,6}\subset\PP(2,2,2,3,3)$. However, the two branch curves will not intersect transversally at $P$ by construction and so this contradicts the hypothesis that $T$ is quasismooth.

Now to fill in the gap we left in the proof that $s_i=t_i=0$, suppose $\delta_1=0$ so that $\al_2=\be_1^{-1}$. Then if $\delta_2=0$ or $\delta_3=0$ this implies $\Delta=0$ which was discounted above. Hence we assume that $\delta_2\delta_3\ne0$ and studying the first and last rows of $C$ we see that this forces $s_5=t_4=0$. However, the remaining two rows of $C$ reduce to $s_4=\al_2t_5$, which no longer has a unique solution!

As a result we get an extension of $\fie$ to \[\widetilde\Fie\colon\PP^5\to\PP(1^4,2^4,3^2)\] where the extra coordinate of weight $2$ is $s_4$ (or equivalently $t_5$). Moreover, the kernel of $\widetilde\Fie^*$ contains equations
\begin{align*}z_1^2-y_1f^2&\in R[s_4]+R[s_4]z_1+R[s_4]z_2\\z_2^2-y_3g^2&\in R[s_4]+R[s_4]z_1+R[s_4]z_2,\end{align*} but the analogue of Corollary~\ref{cor!extend} does not hold unless we insist that $s_4\equiv0$, so that we recover our original hypothesis.

Thus for those particular configurations of degenerate branch curves on $T^\prime_{6,6}\subset\PP(2,2,2,3,3)$ with $\delta_1=0$, there is an extension to some Fano $7$-fold \[V^\prime_{6,6}\subset\PP(1^4,2^4,3^2).\] This does not invalidate the Main Theorem~\ref{thm!nonhell}, since we are looking for Fano $6$-folds $W^\prime_{6,6}\subset\PP(1^4,2^3,3^2)$. However, this is a curious extra stratum of extensions of the K3 surface which merits further investigation.

\subsection{Proof of Corollary~\ref{cor!extend}}\label{sec!pfcor!extend}
To prove the corollary we must calculate the equations extending (\ref{eqn!fieimage2}) and multiples of (\ref{eqn!conic}). First note that \[z_1z_2-fg(y_2+ad-bc)=-fg(bu+cv),\] so if we can write $fg(bu+cv)$ as an expression in the module $R+Rz_1+Rz_2$ then we are done. We must find some $\nu$ in $R\oplus 4R(-3)\oplus R(-4)$ such that \[fg(bu+cv)=\begin{pmatrix}1,&u,&v\end{pmatrix}B\nu.\] Indeed, we can choose the vector $\nu$ such that $\nu_2=bg$, $\nu_3=cf$ and the other $\nu_i=0$. Thus the equation extending (\ref{eqn!fieimage2}) is \[z_1z_2=fg(y_2+ad-bc)-bgz_1-cfz_2.\]

The equations extending (\ref{eqn!conic}) are more complicated. First note from the definition of the matrix $A$ of (\ref{matrix!A}) that \[L_3+L_2u+L_1v=0.\]Thus to write down an equation for $y_iL_3$ in the kernel of $\Fie^*$ we seek some $\nu$ in $R\oplus 4R(-3)\oplus R(-4)$ such that 
\[y_iL_2u+y_iL_1v=\begin{pmatrix}1,&u,&v\end{pmatrix}B\nu.\] Since we used the last column of $B$ to calculate the residual part of $y_iL_3$, to avoid trivial solutions we only use the first $5$ columns of $B$. As previously, the important part is the submatrix \[B^\prime=\begin{pmatrix}f&0&-y_2&y_3\\0&g&y_1&-y_2\end{pmatrix}.\]

Let us calculate the equation for $y_1L_3$. We construct the preimages of $y_1L_2u$ and $y_1L_1v$ under $B$ separately and then sum these two expressions to get the preimage of the residual part. The idea is to try to write down two separate expressions for $y_1y_i$ in terms of $y_if$ and in terms of $y_ig$. With this in mind, consider the resultant matrix 
\begin{equation}\label{matrix!resultant}T=\renewcommand{\arraystretch}{1.2}\left(\begin{array}{ccc|ccc}
1&\al_1&\al_2&&&\\
&1&\al_1&\al_2&&\\
&&1&\al_1&\al_2&\\
\hline
&\be_1&\be_2&1&&\\
&&\be_1&\be_2&1&\\
&&&\be_1&\be_2&1
\end{array}\right).\end{equation}
The matrix $T$ and its inverse have block form \[T=\begin{pmatrix}V_1&V_2\\W_1&W_2\end{pmatrix},\quad T^{-1}=\begin{pmatrix}v_1&v_2\\w_1&w_2\end{pmatrix},\] so that in particular, \begin{align}\label{eqn!identitytop}v_1V_1+v_2W_1&=I_3&v_1V_2+v_2W_2&=0\\
  \label{eqn!identitybottom}w_1V_1+w_2W_1&=0&w_1V_2+w_2W_2&=I_3\end{align} The reason for writing $T$ in block form is that
\[\begin{pmatrix}V_1&V_2\end{pmatrix}\begin{pmatrix}y_1^2\\y_1y_2\\y_1y_3\\y_2y_3\\y_3^2\\ *\end{pmatrix}=\begin{pmatrix}y_1f\\y_2f-\al_1(y_2^2-y_1y_3)\\y_3f\end{pmatrix},\]
where here and elsewhere a star means that entry is irrelevant because it is multiplied by zero. Now if we try to ``invert'' this matrix equation we get an expression for $y_1y_i$ in terms of $y_if$ after a small correction. Multiplying both sides by block $v_1$ and using identities (\ref{eqn!identitytop}) we get
\[\begin{pmatrix}y_1^2\\y_1y_2\\y_1y_3\end{pmatrix}=v_1\begin{pmatrix}y_1f\\y_2f-\al_1(y_2^2-y_1y_3)\\y_3f\end{pmatrix}+v_2W_1\begin{pmatrix}*\\y_1y_2\\y_1y_3\end{pmatrix}+v_2W_2\begin{pmatrix}y_2y_3\\y_3^2\\ *\end{pmatrix}.\] A similar treatment multiplying the bottom half of $T$ by $v_2$ leads to the matrix equation\[\begin{pmatrix}0\\0\\0\end{pmatrix}=v_2\begin{pmatrix}y_1g\\y_2g-\be_2(y_2^2-y_1y_3)\\y_3g\end{pmatrix}-v_2W_1\begin{pmatrix}*\\y_1^2\\y_1y_2\end{pmatrix}-v_2W_2\begin{pmatrix}y_1y_3\\y_2y_3\\ *\end{pmatrix}.\]
Now we can write these two equations in terms of the columns of $B^\prime$ by collecting the terms together appropriately to obtain
\[\setlength\arraycolsep{2pt}\begin{array}{clcll}
y_1Y&=v_1Yf&&+(v_2X_4 + Z_4)(-y_2)&+(v_2X_5 + Z_5)y_3\\
0&=&v_2Yg&+(v_2X_4+Z_4)y_1&+(v_2X_5+Z_5)(-y_2),
\end{array}\]
where \[\quad X_4=\begin{pmatrix}-\be_1y_1-y_3\\-\be_2y_3\\-\be_1y_3 \end{pmatrix},\quad X_5=\begin{pmatrix}\be_2y_1\\\be_1y_1+y_3\\\be_2y_3\end{pmatrix},\quad Y=\begin{pmatrix}y_1\\y_2\\y_3\end{pmatrix},\]
\[Z_4=v_1\begin{pmatrix}0\\\al_1y_2\\0\end{pmatrix}+v_2\begin{pmatrix}0\\\be_2y_3\\0\end{pmatrix},\quad Z_5=v_1\begin{pmatrix}0\\\al_1y_1\\0\end{pmatrix}+v_2\begin{pmatrix}0\\\be_2y_2\\0\end{pmatrix}.\]
The matrices $X_4$, $X_5$ express the terms multiplying $W_1$, $W_2$ above in terms of the columns of $B^\prime$. Similarly $Z_4$ and $Z_5$ express the correction terms involving $y_2^2-y_1y_3$.
Thus multiplying on the left by the matrix $\Lambda_2=\begin{pmatrix}a,&b,&c\end{pmatrix}$ we get \[y_1L_2u=\begin{pmatrix}1,&u,&v\end{pmatrix}B\nu,\] where
\[\nu_2=\Lambda_2v_1Y,\quad \nu_3=\Lambda_2v_2Y,\quad \nu_4=\Lambda_2(v_2X_4+Z_4),\quad \nu_5=\Lambda_2(v_2X_5+Z_5)\]
and \[\nu_1=-\nu_4L_1-\nu_5L_2\] is chosen to cancel extra terms arising from the first row of $B$.

We perform a similar calculation to get an expression for $y_1L_1v$ in the image of $B$. However, this time it is necessary to alter $T$. Let $\sigma$ be the cyclic permutation $(3,4,5,6,1,2)$ of order $3$ acting on the columns of $T$, and let $\sigma^{-1}$ act on the rows of $T^{-1}$. We write $\sigma(T)$ and $\sigma^{-1}(T^{-1})$ in block form as
\[\sigma(T)=\begin{pmatrix}\Vhat_1&\Vhat_2\\\What_1&\What_2\end{pmatrix},
\quad\sigma^{-1}(T^{-1})=\begin{pmatrix}\vhat_1&\vhat_2\\\what_1&\what_2\end{pmatrix}.\]
Then
\[\begin{pmatrix}\What_1&\What_2\end{pmatrix}\begin{pmatrix}y_2y_3\\y_3^2\\ *\\y_1^2\\y_1y_2\\y_1y_3\end{pmatrix}=\begin{pmatrix}y_1g\\y_2g-\be_2(y_2^2-y_1y_3)\\y_3g\end{pmatrix},\]
so that multiplying by $\what_2$ and using permuted versions of identities (\ref{eqn!identitybottom}) we obtain
\[y_1Y=\what_2\begin{pmatrix}y_1g\\y_2g-\be_2(y_2^2-y_1y_3)\\y_3g\end{pmatrix}+\what_1\Vhat_1\begin{pmatrix}y_2y_3\\ *\\ *\end{pmatrix}+\what_1\Vhat_2\begin{pmatrix}y_1^2\\y_1y_2\\y_1y_3\end{pmatrix}.\]
A similar equation arises when we multiply the top half of $\sigma(T)$ by $\what_1$:
\[0=\what_1\begin{pmatrix}y_1f\\y_2f-\al_1(y_2^2-y_1y_3)\\y_3f\end{pmatrix}-\what_1\Vhat_1\begin{pmatrix}y_3^2\\ *\\ *\end{pmatrix}-\what_1\Vhat_2\begin{pmatrix}y_1y_2\\y_1y_3\\y_2y_3\end{pmatrix}.\]
Then separate out these two equations as expressions in the columns of $B$
\[\setlength\arraycolsep{2pt}\begin{array}{clcll}
0&=\what_1Yf&&+(\what_1\widehat X_4+Z_4)(-y_2)&+(\what_1\widehat X_5+Z_5)y_3\\
y_1Y&=&\what_2Yg&+(\what_1\widehat X_4+Z_4)y_1&+(\what_1\widehat X_5+Z_5)(-y_2),\end{array}\]where
\[\widehat X_4=\begin{pmatrix}\al_1y_1\\y_1+\al_2y_3\\\al_1y_3\end{pmatrix},\quad \widehat X_5=\begin{pmatrix}-\al_2y_1\\-\al_1y_1\\-y_1-\al_2y_3\end{pmatrix}\] and $Z_4$, $Z_5$ are as above. We multiply on the left by $\Lambda_1:=\begin{pmatrix}b,&c,&d\end{pmatrix}$ to obtain an expression for $y_1L_1v$ in the image of $B$. The preimage $\nuhat$ is the vector
\[\nuhat_1=-\nuhat_4L_1-\nuhat_5L_2,\quad\nuhat_2=\Lambda_1\what_1Y,\quad\nuhat_3=\Lambda_1\what_2Y,\]
\[\nuhat_4=\Lambda_1(\what_1\widehat X_4+Z_4),\quad
\nuhat_5=\Lambda_1(\what_1\widehat X_5+Z_5),\]
Hence \[y_1(L_2u+L_1v)=\begin{pmatrix}1,&u,&v\end{pmatrix}B(\nu+\nuhat)\] is the residual part to $y_1L_3$ and so we can write out an equation in $R+Rz_1+Rz_2$:
\[y_1L_3+(\nu_1+\nuhat_1)+(\nu_2+\nuhat_2)z_1+(\nu_3+\nuhat_3)z_2=0.\]

The calculation of $y_2L_3$, $y_3L_3$ requires further cyclic permutations of the columns of $T$. We do not write out the details here but it follows the same pattern as the calculations above.

\section{Surfaces with $p_g=1$ and $K^2=2$}\label{sec!surfaces}
This brief final section consists of the following application of our $6$-fold extensions to surfaces of general type.
\begin{thm}\label{thm!keyvar}There is a $16$ parameter family of surfaces $Y$ of general type with $p_g=1$, $q=0$, $K^2=2$ and no torsion, each of which is a complete intersection of type $(1,1,1,2)$ in a Fano $6$-fold $W\subset\PP(1^4,2^4,3^4)$ with $10\times\frac12$ points.
\end{thm}
\begin{pf} To obtain $Y$ from $W$ take $3$ transverse hyperplane sections of weight $1$ and one hypersurface section of weight $2$, avoiding the isolated orbifold $\frac12$ points. Since $W$ is quasismooth, by the adjunction formula the surface $Y$ has $\omega_Y=\Oh_Y(1)$ and is smooth. Furthermore it is clear from the construction of $Y$ that \[p_g(Y)=h^0(Y,K_Y)=h^0(W,-K_W)-3=1.\] Consider $Y$ as a quadric section of a Fano $3$-fold $W^3$. Then the standard short exact sequence \[0\to\Oh_{W^3}(-2)\to\Oh_{W^3}\to\Oh_Y\to0,\] implies that $H^1(\Oh_Y)=0$ by Kodaira vanishing, so $Y$ is regular. Finally the Riemann--Roch formula gives $K_Y^2=2$.

Theorem~\ref{thm!nonhell} says that the family of Fano $6$-folds depends on the same number of moduli as the family of symmetric determinantal quartics, which is $9$. Furthermore, naively counting the number of choices for linear and quadric sections of $W$ suggests that we have a $9+3+4=16$ parameter family of surfaces $Y$.
\end{pf}

This agrees with the expected dimension of the moduli space of surfaces with $p_g=1$, $K^2=2$, which suggests that we have constructed the general surface (see \cite{CDE} for further justification). However, we have not proved that \emph{every} such surface is contained in a Fano $6$-fold $W$ as a weighted hyperplane section, only that the canonical curve section $D$ is.

\begin{small}\end{small}
\bigskip\noindent
Stephen Coughlan, \\
Math. Dept., Sogang University, \\
Sinsu-dong, Mapo-gu, \\
Seoul, Korea \\
e-mail: stephencoughlan21@gmail.com \\

\end{document}